\documentclass{amsart}

\usepackage{amsmath}
\usepackage{amssymb}
\usepackage[all]{xy}

\newcommand{\D}{\mathbb{D}}
\newcommand{\C}{\mathbb{C}}

\newcommand{\PP}{\mathbb{P}}

\newcommand{\N}{\mathbb{N}}
\newcommand{\Z}{\mathbb{Z}}

\newcommand{\lgra}{\longrightarrow}
\newcommand{\ra}{\rightarrow}
\newcommand{\dsps}{\displaystyle}

\newcommand{\exc}{\setminus}

\newcommand{\td}{\widetilde}

\def \pc{almost-complex}

\def \hq{hyperbolic}
\def \hqe{hyperbolic }
\def \hle{holomorphic }

\def \hcte{hyperbolicity }

\def \hl{holomorphic}

\def \ste{structure }

\def \spe{symplectic }

\def \pp{\partial}
\def \f{\frac}
\def \dd{{\rm d}}

\def \mcl{\mathcal}
\def \be{\begin{equation}}
\def \ee{\end{equation}}
\def \beq*{\begin{equation*}}
\def \eeq*{\end{equation*}}
\def \ba{\begin{eqnarray}}
\def \ea{\end{eqnarray}}
\def \ba*{\begin{eqnarray*}}
\def \ea*{\end{eqnarray*}}

\newtheorem{prop}{Proposition}[section]

\newtheorem{cor}{Corollary}[section]

\newtheorem{defi-thm}{Definition-Theorem}[section]

\newtheorem{theorem}{Theorem}[section]
\newtheorem{lemma}[theorem]{Lemma}

\theoremstyle{definition}
\newtheorem{definition}[theorem]{Definition}

\theoremstyle{remark}

\numberwithin{equation}{section}

\begin{document}

\title{Hyperbolicity of holomorphic foliations with parabolic leaves}

\author{Anne-Laure Biolley}
\address{Department of Mathematics, University of Toronto, Ontario M5S
3G3}

\email{alb@math.toronto.edu}
\thanks{Work supported by Ecole Polytechnique (France) and
University of Toronto (Canada).}

%    General info

\date{July 30, 2004.}

\dedicatory{I would like to thank Etienne Ghys for sharing with me
his ideas on this question and planting the seeds of this paper in
me, and my PhD supervisor, Claude Viterbo, for his remarks and
advice.}


\begin{abstract}
In this paper, we study a notion of hyperbolicity for \hle
foliations with $1$-dimensional parabolic leaves, namely the
non-existence of \hle cylinders along the foliation -- \hle maps
from $\D^{n-1} \times \C$ to the manifold sending each $\{*\}
\times \C$ to a leaf. We construct a tensor, that is a leaf-wise
\hle section of a bundle above the foliated manifold, which
vanishes on an open saturated set if and only if there exists such
a cylinder through each point of this open set. Thanks to this
study, we prove that open sets saturated by compact leaves are
union of \hle cylinders along the foliation. We also give some
more specific results and examples in the case of a compact
manifold, and of foliations with one compact leaf or without any
compact leaves.\end{abstract}

\maketitle
\section*{Introduction}
This paper is a first step in the study of hyperbolicity of
foliations. In the general case of complex manifolds, a way of
tackling complex hyperbolicity is to study the existence of \hle
curves from $\C$ to the manifold: if the only such curves are the
constant-curves, then the manifold is said to be Brody-\hqe (see
notably \cite{lang}, \cite{alb} or \cite{alb1}). However, in the
special case of foliations, as soon as there exists a parabolic
leaf (\emph{i.e.} a leaf whose universal covering is the complex
plane $\C$) this question becomes trivial. That is why in the case
of foliations with parabolic leaves, we need to study some more
appropriate and relevant notions of hyperbolicity. More precisely,
this paper studies the existence of \hle maps from $\D^{n-1}
\times \C$ to the manifold along the foliation, \emph{i.e.} maps
which send each $\{*\} \times \C$ to a leaf ($\D^{n-1}$ denoting
the ball of radius $1$ in $\C^{n-1}$). Then, if there does not
exist any such map which is immersed (there always exist the
trivial ones taking values into a single leaf), the foliated
manifold is called \hq.\\
\par
This question is linked with the issue of measure-hyperbolicity
(see \cite{lang} for the definitions of the different notions of
hyperbolicity and their links). As the Kobayashi \hcte corresponds
to the non-degeneracy of the Kobayashi pseudo-distance (see
\cite{lang}), the measure-hyperbolicity corresponds to the
non-degeneracy of the Kobayashi measure. In a similar way, as the
existence of non-constant \hle maps from $\C$ to a complex
manifold implies its non-Kobayashi hyperbolicity, the existence of
non-trivial \hle maps from $\D^{n-1} \times \C$ implies its
non-measure-hyperbolicity. This is one of the motivations behind
the study of the existence of such maps.\\
\par
For this study, we will consider foliations by parabolic leaves
(so that $\{*\} \times \C$ could be sent onto a leaf).\\

%Thanks to the analysis of connections defined on the tangent bundle
%to the foliation, we defined an invariant of the foliation, which
%vanishes if and only if such a map exists. This provides us with
%conditions in terms of existence of nonzero \hle sections of a
%bundle \hle along the foliation. Moreover, this allows to get some
%results in the theory of foliations by Riemann surfaces.
\par
So let us consider a complex manifold $X$ of dimension $n$ (real
dimension $2n$) provided with a $1$-dimensional \hle foliation
$\mathcal{F}$ whose leaves are all parabolic (the leaf passing
through a point $x$ will be denoted by $L_x$). We would like to
find conditions under which there exists a \hle cylinder along the
foliation:
\begin{definition}
A \emph{\hle cylinder along the foliation} is an immersed \hle map
$F~: \D^{n-1} \times \C \ra X$ such that, for any $x \in \D$,
$\{x\}\times \C$ is sent onto a leaf of the foliation
\end{definition}

%(The \emph{immersion}-condition will be implicit in the
%following).

We construct in this paper an invariant of the foliation: a
leaf-wise \hle section of a \hle bundle of the foliation, which
vanishes if and only if there exist such \hle cylinders along the
foliation. More precisely:
\begin{theorem}
There exists a tensor $\Gamma$ which is an invariant associated to
the foliation:
\begin{equation*}
\Gamma_{|x}: T_{x}X \times T_{x}\mathcal{F} \ra \C, \ \forall x
\in X
\end{equation*}
and which vanishes on an open saturated set $U$ if and only if
through any point of $U$ there exists a \hle cylinder $F: \D^{n-1}
\times \C \ra U$ along the foliation. \\
It is $\C$-antilinear with respect to the first variable and
$\C$-linear with respect to the second one.\\
 Moreover, $\Gamma$ is a section of the bundle
$E=\Lambda^{(0,1)} T^*\mathcal{N} \otimes \Lambda^{(1,0)}
T^*\mathcal{F}$ over $(X, \mathcal{F})$, which is \hle along the
leaves.
\end{theorem}
\noindent In this theorem, and more generally in this paper,
$T\mathcal{N}$ denotes the transverse (or normal) vector bundle of
the foliation: $T \mathcal{N}=TX/T\mathcal{F}$. The usual
vocabulary and properties of \hle foliations are recalled in
Appendix \ref{rappel}.\\
\par
This analysis allows us to get some new results in the theory of
foliations. Notably, in the case of a foliation with compact
leaves, we can deduce that through each point, there exists a \hle
cylinder along the foliation (result which was already known in
the case of Stein manifolds foliated by parabolic leaves). This
implies that, in this case, the leaves depend holomorphically on
the transversal disk (see
section \ref{compactleaves} for more details).\\

In order to simplify notations, the analysis of this paper will
first be conducted in the case of complex dimension $2$. But, as
it will be underlined, section \ref{higherdim}, the results can
straightforwardly be generalized to the case of higher dimensions.
Let us detail the structure of
this paper.\\

First, in section \ref{constructF}, we introduce a particular set
$\mathcal{C}$ of maps $F: \D \times \C \ra X$ which naturally
appear as candidates for this analysis. And we prove that if a map
of the desired type exists, then it necessarily belongs to this
set $\mathcal{C}$. Thus, the problem reduces to the question of
the holomorphy of the maps in $\mathcal{C}$.\\
\par
To tackle this problem, we associate to any map $F \in
\mathcal{C}$, a complex-valued function $\omega$ on $\D \times
\C$, which is \hle along the leaves $\{*\} \times \C$, and
measures the lack of holomorphy of $F$: this function is zero if
and only if $F$ is \hl. This function $\omega$ depends on the map
$F$. However, we note that if one of the maps in $\mathcal{C}$ is
\hl, then so are all the other maps in $\mathcal{C}$ (whose images
are included in the one of $F$). This leads us to look for an
invariant on $X$ measuring the lack of holomorphy of any maps in
$\mathcal{C}$.\\
\par
 To deal with this, we study how $\omega$, and the
local differential forms on $X$ naturally associated to it, depend
on $F$ (section \ref{changt}, section \ref{paragdef} and
\ref{paragchang} respectively).\\
\par
This leads us to construct such an invariant $\Gamma$, a tensor on
$X$, independent of any choice of $F$, which is zero on an open
saturated set if and only if any map $F \in \mathcal{C}$ taking
values in this open set is \hl. This invariant is described in
section \ref{invt}. It is an invariant associated to the foliation
which detects the existence of \hle cylinders along the foliation.
More precisely, let us state our result:
\begin{theorem}
To any $F\in \mathcal{C}$ is associated a tensor, defined on an
open set $V_F$ around $x_0=F(0,0)$,
\begin{equation*}
{\Gamma_F}_{|x}: T_x X \times T_x \mathcal{F} \ra \C, \ \forall x
\in V_F.
\end{equation*}
This tensor does not depend on $F$, \emph{i.e.} if $F' \in
\mathcal C$ and $V_F \cap V_{F'} \neq \emptyset$, then $\Gamma_F=
\Gamma_{F'}$ on $V_F \cap
V_{F'}$.\\
Thus, a tensor can be defined on $X$ by $\Gamma=\Gamma_F$ on each
open set $V_F$. It is $\C$-antilinear with respect to the first
variable and $\C$-linear with respect to the second one. This is
an invariant of the foliation, that is zero if and only if the
maps in $\mathcal{C}$ are \hl. More precisely, for any $F \in
\mcl{C}$, $\Gamma$ is zero on
$F(\D\times \C)$  if and only if $F$ is \hl.\\
Moreover, $\Gamma$ is a section of the line bundle
$E=\Lambda^{(0,1)} T^* \mathcal{N} \otimes \Lambda^{(1,0)} T^*
\mathcal{F}$ over $(X, \mathcal{F})$ and it is \hle along the
leaves.
\end{theorem}

Thus if $\Gamma=0$ on the saturation $U$ of a transverse disk $D$,
then there exists a \hle cylinder along the foliation $F: \D
\times \C \ra X$ that sends $\D \times \{0\}$ biholomorphically on
$D$. In fact, any $F \in \mathcal{C}$ whose image is included in
$U$ is such a \hle cylinder.\\
\par
Moreover, if $\Gamma=0$ on $X$, all the maps in $\mathcal{C}$ are
\hl. This, in particular, is satisfied when the line bundle $E$
over $(X, \mcl{F})$ does not admit non-zero leaf-wise holomorphic
sections, since $\Gamma$ is a leaf-wise holomorphic section of
$E$.\\
\par

Finally, we study some particular cases and examples. First, we
prove that if all the leaves of $\mathcal{F}$ are compact, then
$\Gamma$ identically vanishes, and therefore all the maps $F\in
\mathcal{C}$ are \hle (section \ref{compactleaves}). Similarly, if
$\mathcal{F}$ has one compact leaf $L_0$ with a finite holonomy
group, then $\Gamma$ vanishes on an open saturated set $U$ around
$L_0$, and all the maps $F \in \mcl{C}$ taking values in $U$ are
\hl. In the case where $\mathcal{F}$ has one compact leaf with an
infinite holonomy group, we study the action of the holonomy group
on the complex-valued function $\omega$. We also prove that in the
case where the manifold $X$ is compact, even though the leaves are
not, if the bundle $E$ has a negative curvature along the leaves
then the invariant $\Gamma$ vanishes. We also give an example of a
manifold with non-compact leaves whose invariant $\Gamma$ is
non-zero and therefore does not admit any \hle cylinder along the
foliation (section \ref{nontrivialex}).

\section{Holomorphy of cylinders along the foliation}
\subsection{Cylinders along the foliation}
\label{constructF} Let us consider $X$ a $2$-dimensional complex
manifold provided with a \hle foliation of complex dimension $1$
whose leaves are all parabolic. We aim to study the existence of
\hle maps from $\D \times \C$ to $X$ along the foliation.\\

Let us first note that if such a \hle map $F:\D\times \C \ra X$
along the foliation exists, then $D:=F(\D \times \{0\})$ is a \hle
disk transverse to the foliation. Moreover, for every $x \in \D$,
$t_x:=\dd F_{(x,0)} \begin{pmatrix} 0 \\ 1 \end{pmatrix} \in T
\mcl{F}$, and $x \ra t_x$ is a \hle trivialization of
$T\mathcal{F}$ along $D$.\\

Thus the cylinder-maps $F:\D \times \C \ra X$ along the foliation
(\emph{i.e.} a map sending $\{ *\} \times \C$ onto a leaf) such
that $D=F(\D \times \{0\})$ is a \hle disk transverse to the
foliation and $t_x=\dd F_{(x,0)}
\begin{pmatrix} 0
\\ 1 \end{pmatrix}$ is a \hle trivialization of
$T\mathcal{F}$ along $D$, appear as the natural and only
candidates for our study. The question is: are these maps \hl?\\

We will denote by $\mcl{C}$ the sets of all these maps. Let us
first describe the one-and-one correspondence between the data of
such a map $F$, and the data of the transverse disk $D$ and of the
trivialization $t$. \\

Let $D$ be a \hle disk transverse to the foliation, parametrized
\hl ally by the unit disk $\D$ with which it will be identified.
Along $D$, we fix a trivialization $t:\D \ra T\mathcal{F}_{|D} $
of the tangent bundle $T\mathcal{F}$ of the foliation. For each $x
\in \D$, the leaf $L_x$ passing through $x$ is parabolic. Thus, if
$\td{L_x}$ denotes its universal cover, $\td{L_x} \simeq \C$. Then
for any fixed covering map  $\psi_x~:\td{L_x} \ra L_x$, there
exists an unique bi-holomorphism $F^0_x~: \C \ra \td{L_x}$ such
that $F_x:=\psi_x \circ F^0_x$ satisfies $F_x(0)=x$ and $\dd
F_0(1)=t_x$. Thus we get a surjective immersion

\begin{equation*}
F~: \begin{cases} \D \times \C & \lgra X \\
                 (x,y) & \lgra F_x (y)
\end{cases}
\end{equation*}
which is by construction holomorphic with respect to the second
variable $y$ and holomorphic on $\D \times \{0\}$.\\

It is the only map from $\D \times \C$ along the foliation such
that $F(\D \times \{ 0\})=D$ (the parametrization of $D$ via $\D$
being fixed) and $\dd F_{|\D \times \{0\}} \begin{pmatrix} 0
\\ 1 \end{pmatrix} =t$.\\
\par
If one of these maps in $\mcl{C}$ is \hl, then it is a map of the
desired kind: a \hle cylinder along the leaves.\\

Reciprocally, if such a \hle map $F$ along the foliation exists,
then, as noted in the beginning of this section, this map $F$
belongs to $\mathcal{C}$. That is why, for the analysis of our
problem, we only have to study the holomorphy of these maps $F \in
\mathcal{C}$.

\subsection{Holomorphy of the maps $F\in \mathcal{C}$}

Let $F$ be a map in $\mathcal{C}$, and let us consider $J'=F^* J$,
with $J$ the complex structure on $X$. This \ste $J'$ can be
decomposed as $J'= J_0 + \Omega$ with $J_0$ the standard complex
structure on $\D \times \C$ and $\Omega$ a tensor in $End(T(\D
\times \C))$. This tensor $\Omega$ measures the lack of holomorphy
of $F$, since it is zero if and only if $J'=J_0$, which means if
and only if $F$ is \hl.\\

This tensor $\Omega$ takes values in $T \mcl{F}_0$ because $J'$
and $J_0$ coincide on $\C \times \{ 0 \} \subset T(\D \times \C)$
(since $F$ has been constructed along the \hle transverse disk
$D$; for more details see the Appendix \ref{projhol}). Moreover,
it vanishes on vectors tangent to the foliation (because $F$ is by
construction \hle with respect to $y$).\\

Furthermore, $(J')^2=-Id$, which implies that $\Omega J_0 = J_0
\Omega$ and $\Omega$ is antilinear. Thus, there exists a function
$\omega~: \D \times \C \ra \C$ with $\omega_{|\D
  \times \{0\}} = 0$ such that

\begin{equation*}
\Omega_{(x,y)} \begin{pmatrix} v_x \\ v_y
\end{pmatrix}=\begin{pmatrix} 0 \\ \omega(x,y) \overline{v_x} \end{pmatrix}.
\end{equation*}

\par
This function $\omega$ satisfies:

\begin{lemma}
The function $\omega : \D \times \C \ra \C $  defined above
vanishes if and only if $F$ is holomorphic. It is \hle with
respect to the second variable, that is to say for all $x \in \D$,
the functions $\omega(x,.)$ are holomorphic on $\C$. Moreover
these functions vanish on 0 and their derivatives also vanish on
0.
\end{lemma}
\begin{proof}
First, by definition of $\Omega$, $F$ is \hle if and only if
$\Omega=0$, or equivalently if and only if $\omega=0$.\\

\noindent As $F$ is holomorphic along $\D \times\{0\}$, for all $x
\in \D$,
$\omega(x,0)=0$.\\

\noindent The leaf-wise holomorphy of $\omega$ comes from the
integrability of the \ste $J'$. Indeed this property translates
as: the Nijenhuis vector field $N_{J'}(Z,Y)$ vanishes for all
vector fields $Z$ and $Y$ (or, equivalently, for any constant
vector fields $Z$ and $Y$). Since $J_0$ is integrable, and $J'
=J_0 +\Omega$,
\begin{eqnarray*}
N_{J'}(Z,Y) &=& [\Omega Z, J_0 Y]+[J_0 Z, \Omega Y] + [\Omega Z,
\Omega Y] -\Omega [J_0 Z, Y] - \Omega [\Omega Z, Y] \nonumber \\ &
& - J_0 [\Omega Z, Y] - \Omega [Z, J_0 Y] - J_0 [Z, \Omega Y] -
\Omega [Z, \Omega Y].
\end{eqnarray*}
Let us take $Z=(z_1, z_2)$ and $Y=(y_1,y_2)$ constant on $\D
\times \C$. Since $\nabla J_0=0$, and since $\Omega$ is zero on
$\D \times \{0\}$ and takes values in $\{0\} \times \C$,
$N_J'(Z,Y)=0$ is equivalent to :
\begin{eqnarray*}
d\omega(J_0 Y) \bar{z_1} - d\omega(J_0 Z) \bar{y_1} +d\omega(0,
\omega \bar{y_1}) \bar{z_1} & \nonumber \\
-d\omega(0, \omega \bar{z_1}) \bar{y_1} -J_0 d\omega(Y) \bar{z_1}
+ J_0 d\omega(Z) \bar{y_1} & =0.
\end{eqnarray*}
Firstly fixing $y_1=0$ and $x_1 \neq 0$, we notice that if $J'$ is
integrable then $\omega$ has to be holomorphic with respect to
the second variable.\\
Reciprocally, if $\omega$ is holomorphic with respect to the
second variable, one can straightforwardly check that
$N_{J'}(Z,Y)=0$ for all constant vector fields $Z$ and $Y$, and
therefore for every vector field. Thus, $J'$ (or equivalently $J$)
is integrable if and only if $\omega$ is holomorphic with respect
to the second variable.\\

\noindent Let us now translate the holomorphy of the map  $x \ra
t_x = dF_{(x,0)}
\begin{pmatrix} 0 \\ 1 \end{pmatrix}$.
Let us consider a distinguished open set $U$ containing $F(\D
\times \{ 0\})$ on which there exists a \hle foliated chart:
$\psi~: U \ra \D \times \D$ $(J,J_0)$-\hle (with $J_0$ the
standard structure on $\D \times \D$), sending $D$ on $\D \times
\{0\}$. (The definitions of the involved notions and terms are
recalled Appendix \ref{rappel}).

\noindent Then let us define $\alpha:=\psi \circ F$ on
$F^{-1}(U)$. \xymatrix{ F^{-1}(U) \ar@{^{(}->}[d] \ar[r]^{F} & U
\ar@{^{(}->}[d] \ar[r]^{\psi} &
\D \times \D \\
\D \times \C \ar[r]^F & X & }

\noindent This map $\alpha$ can be written as
$\alpha(x,y)=(\alpha_1(x), \alpha_2(x,y))$, with $\alpha_1$ a
bi-holomorphism of $\D$ and $\alpha_2$ a map \hle with respect to
$y$. Moreover, by writing that $\alpha$ is holomorphic for the
complex structures $J'$ and $J_0$,

\beq* \dd \alpha \left( (\Omega +J_0) \begin{pmatrix} v_x \\ v_y
\end{pmatrix}\right) = J_0 \dd \alpha \begin{pmatrix} v_x \\ v_y
\end{pmatrix}
\eeq* \noindent we get \be \label{omegalapharel}2 i \frac{\partial
\alpha_2}{\partial
  \bar{x}}=\omega \frac{\partial \alpha_2}{\partial y}.\ee
  So, by differentiating:
\begin{equation}
\label{condalpha}
2 i \frac{\partial^2 \alpha_2}{\partial y \partial
  \bar{x}}=\frac{\partial \omega}{\partial y}
\, \frac{\partial \alpha_2}{\partial y} + \omega\frac{\partial^2
  \alpha_2}{\partial y^2} .
\end{equation}

\noindent Besides, $x \ra t_x$ is \hl, that is to say, in the
trivialization, $x \ra \frac{\partial  \alpha_2}{\partial y}(x,0)$
is \hl, which reads as $\frac{\partial^2 \alpha_2}{\partial y
\partial \bar{x}}(x,0)=0$. And since $\omega(x,0)=0$, we
deduce from this that $\frac{\partial
  \omega}{\partial y}=0$ on $\D \times \{ 0\}$.\\

\noindent Summing up, $\omega(x,.)$ are \hle functions on $\C$
which vanish on 0 and whose derivatives also vanish on 0.
Furthermore $F$ is \hle if and only if these functions are
identically zero, or equivalently if and only if $\frac{\partial
  \omega}{\partial y}=0$ (because $\omega=0$ on $\D \times \{0\}$ and
it is \hle with respect to $y$).
\end{proof}

So to summarize, the function $F$ is \hle if and only if the
associated function $\omega$ vanishes, or equivalently, if and
only if $\frac{\partial  \omega}{\partial y}=0$, or also if and
only if $\frac{\partial^2  \omega}{\partial y^2}=0$. We are now
going to study more this function $\omega$.

\subsection{How $\omega$ depends on $F$}
\label{changt} The goal of this section is to see how the function
$\omega$ varies when one considers another function in
$\mcl{C}$.\\

So let us consider another $\tilde{F} \in \mathcal{C}$ constructed
as in \ref{constructF}, from the data of a \hle disk $\tilde{D}$
transverse to the foliation, and of a \hle trivialization
$\tilde{t}$ of $T\mathcal{F}$ along $\tilde{D}$. Then as
previously, we can associate to it the tensor $\tilde{\Omega}$ and
the function $\tilde{\omega}$ measuring the lack of holomorphy of
$\td{F}$. We would like to compare the function $\omega$
associated with $F$, with this function $\td{\omega}$. More
precisely, if $F(x,y)= \td{F}(\td{x}, \td{y})$, we'd like to
compare $\omega(x,y)$ with $\td{\omega}(\td{x}, \td{y})$.\\
\par
Since $D$ and $\td{D}$ are two transverse \hle disks, a classical
construction explained in the Appendix (lemma
\ref{translationalongleaves}) provides us with a bi-holomorphism
$\theta_1 : D_1 \ra D_2$, with $D_1 \subset \D$ and $D_2 \subset
\D$, such that for any $x$, $L_{\tilde{F}(x,0)}=
L_{F(\theta_1(x),0)}$. So, without any loss of generality (since
we are only interested in the intersection of the images of $F$
and $\td{F}$), possibly restricting the two transverse \hle disks
$D$ and $\td{D}$, we can consider $\theta_1$ as a bi-holomorphism
of $\D$ such that for any $x\in \D$, $\theta_1(x) \in \D$
satisfies $L_{\tilde{F}(x,0)}=
L_{F(\theta_1(x),0)}$.\\

Then, for any $x \in \D$, we can  consider the bi-holomorphism of
$\C$ : $\theta_x=F_{\theta_1(x)}^{-1}\circ \tilde{F}_x$ (actually
$({F^0}_{\theta_1(x)})^{-1}\circ \tilde{F^0_x}$). There exist some
complex numbers $a(x)$, and $b(x)$ such that it can be
written: $\theta_x(y)=a(x)\, y + b(x)$.\\

And we define $\theta~:\D \times \C \ra \D \times \C$ (locally
$\theta = F^{-1} \circ \tilde{F}$) by
 $\theta(x,y) =(\theta_1(x), \theta_2(x,y))$ with $\theta_2(x,y)
=\theta_x(y)$. Thus we prove:
\begin{lemma}
\label{deftheta}
 There exists a diffeomorphism $\theta: \D \times
\C \ra \D \times \C$ satisfying:
\begin{itemize}
\item $\td{F}=F \circ \theta$ \hspace*{2cm} \xymatrix{ \D \times
\C \ar[r]^{\td{F}} \ar[d]_{\theta} & X \\ \D \times \C \ar[ru]_{F}
& }

\item $\theta(x,y) =(\theta_1(x), \theta_2(x,y))$, with $\theta_1$
a bi-holomorphism of the disk $\D$, and for any $x \in \D$
$\theta_2(x,.)$ a bi-holomorphism  of $\C$.
\end{itemize}
Moreover,
\begin{equation}
\label{unepetiteeq}
\frac{\partial
  \theta_2}{\partial y} \, \tilde{\omega}(x,y)=
\omega(\theta(x,y)) \, \overline{\frac{\partial
  \theta_1}{\partial x}} + 2 i \frac{\partial \theta_2}{\partial
\bar{x}}.
\end{equation}
Thus, the change in the function $\omega$ is measured by the lack
of holomorphy with respect to $x$ of the function $\theta_2$.
\end{lemma}
\begin{proof}
Only (\ref{unepetiteeq}) remains to be proven. Let us compare
$\Omega$ and $ \td{\Omega}$ : $J_0+\tilde{\Omega} = (\tilde{F})^*
J= (F\circ \theta)^* J =\theta^* (J_0+\Omega)$. Therefore
$\tilde{\Omega}=\theta^* \Omega + \theta^* J_0 - J_0$, so
$\theta_* \circ \td{\Omega}= \Omega \circ \theta_* + J_0 \circ
\theta_* - \theta_* \circ J_0$, which exactly writes as
(\ref{unepetiteeq}).\\
The last remark of the lemma is explained since, possibly by
reparametrizing $\td{D}$ by a bi-holomorphism of $\D$, we can
suppose $\frac{\partial  \theta_1}{\partial x}=1$. Besides the
term $\frac{\partial \theta_2}{\partial y}$ is a non-zero
multiplicative term -- corresponding to a reparametrization via a
bi-holomorphism of $\C$ of the maps $\td{F}_x$ (but not
necessarily \hle with respect to $x$). So the most relevant term
in (\ref{unepetiteeq}) is $\frac{\partial \theta_2}{\partial
\bar{x}}$.
\end{proof}

Thus, if $F$ is not \hl, then the pull-back in $\D \times \C$, via
$F$, of the \hle disk $\td{D}=\td{F}(\D \times \{0\})$ and of the
\hle trivialization $\td{t}=\dd \td{F}_{(x,0)} \begin{pmatrix} 0
\\ 1 \end{pmatrix}$  will not be \hl. Besides, since
locally $\theta=F^{-1} \circ \td{F}$, the pull-back by $F$ of
$\td{D}$ is $\theta(\D \times \{0\})$, \emph{i.e.} the graph of
$b(x)$; and the pull-back by $F$ of $\td{t}$ is $\dd
\theta_{(x,0)}
\begin{pmatrix} 0 \\ 1 \end{pmatrix}$, \emph{i.e} the graph of
$a(x)$. So the functions $a$ and $b$, and so also the function
$\theta_2 = a(x) y+ b(x)$ are not \hle with respect to the
variable $x$. This lemma tells us that, in this case, in fact, the
change in the function $\omega$ is roughly measured by the
non-holomorphy of
$\theta_2$.\\
\par
Moreover, by the same argument, it is interesting to note that if
$F$ is \hl, then $F$ determines (locally) a \hle trivialization of
the foliation. Thus, both the disk $\tilde{D}$ and the
trivialization $\tilde{t}$ read holomorphically in this
trivialization. This implies that $a(x)$ and $b(x)$ are \hl, and
so is $\theta_2$. So, in regards of (\ref{unepetiteeq}), $\td{F}$
is \hl. Thus
\begin{lemma} If $F \in \mcl{C}$ is \hl, then any $\tilde{F} \in
\mcl{C}$ (whose image is included in the one of $F$) is \hle too.
\end{lemma}

In order to better understand the effect of changing the function
$F$, let us study what the possible changes are, \emph{i.e} the
set of functions $\theta$ which can be associated to a change of
the function $F$ in another function $\td{F} \in \mcl{C}$. This is
particularly relevant since, as we have seen, this function
$\theta$ determines the change in $\omega$ (see
(\ref{unepetiteeq})). So, let us consider some map
\begin{equation*}
\theta~: \begin{cases} \D \times \C & \lgra \D \times \C \\
(x,y) & \lgra (\theta_1(x), \theta_2(x,y)) \end{cases}
\end{equation*}
such that $\theta_1$ is a bi-holomorphism of $\D$ and, that for
every $x \in \D$ the map $y \ra \theta_2(x,y)$ is a
bi-holomorphism of $\C$ (and so of the form $a(x)y+b(x)$). Then,
what are the necessary and sufficient conditions such that
$\tilde{F}:=F \circ \theta$ belongs to $\mcl{C}$ (\emph{i.e.} that
it might be obtained as in \ref{constructF})?\\

\par
These conditions are:

\begin{enumerate}
\item
\label{cond1}
 $F(\theta(\D \times \{0\}))$ is a \hle disk \emph{i.e.}
$F_*(\theta_* \begin{pmatrix} i \\ 0 \end{pmatrix})=JF_*(\theta_*
\begin{pmatrix} 1 \\ 0 \end{pmatrix})$ on $\D \times \{0\}$, which
is equivalent to $-2 i \frac{\partial \theta_2}{\partial
\bar{x}}(x,0) = \omega(\theta(x,0)) \, \overline{\frac{\partial
\theta_1}{\partial x}(x,0)}$, $\forall x$.

\item $x \ra \dd \tilde{F}_{(x,0)}\begin{pmatrix} 1 \\ 0
\end{pmatrix}$ is \hl, \emph{i.e.} $x
  \ra \dd F_{\theta(x,0)}\begin{pmatrix} 0 \\ 1
  \end{pmatrix}\, \frac{\partial \theta_2}{\partial y} (x,0)$ is \hl.
\label{cond2}
\end{enumerate}

Let us assume (possibly by first considering some intermediate
disks) that $\tilde{D}:=F(\theta(\D \times \{0\}))$ is
sufficiently close to $D$, or more precisely that it is included
in a distinguished open set $U$. Let us fix a foliated chart
$\psi~: U \ra \D \times \D$. How is condition (\ref{cond2})
expressed in this trivialization?\\

Keeping the previous notation  $\alpha=\psi \circ F$ on
$F^{-1}(U)$, condition (\ref{cond2}) reads in this chart: the map
$x \ra \frac{\partial
  \alpha_2}{\partial y}(\theta(x,0))  \, \frac{\partial
  \theta_2}{\partial y}(x,0)$ is \hl. By writing that the
$\bar{\partial}$ of this function must be zero and by using the
equality (\ref{condalpha}) and condition  (\ref{cond1}), we get
that (if condition (\ref{cond1}) is satisfied) condition
(\ref{cond2}) is equivalent to:

\begin{equation*}
\left( \frac{\partial \alpha_2}{\partial y}\circ \theta \right) \,
\left(\frac{\partial
  \omega}{\partial y} \circ \theta \, \overline{\frac{\partial
  \theta_1}{\partial x}} \,  \frac{\partial \theta_2}{\partial y} + 2 i
\frac{\partial^2 \theta_2}{\partial y \partial \bar{x}}\right)=0 \
{\rm on}
  \ \D \times \{0\}.
\end{equation*}

Since $\frac{\partial \alpha_2}{\partial y}$ never vanishes, and
since $\frac{\partial \omega}{\partial y} \circ \theta$ is the
only term depending on $y$ (neither does $\theta_1$ depend on $y$,
nor $\frac{\partial \theta_2}{\partial y}=b(x)$), we get:

\begin{equation}
\label{relation}
\frac{\partial
  \omega}{\partial y} (\theta(x,0))\, \overline{\frac{\partial
  \theta_1}{\partial x}}(x,y)\,  \frac{\partial \theta_2}{\partial y}(x,y) + 2 i
\frac{\partial^2 \theta_2}{\partial y \partial \bar{x}}(x,y)=0 \ {\rm
  on \ } \  F^{-1}(U).
\end{equation}

\par
So, to summarize, if we assume (possibly by composing with a
bi-holomorphism of the disk) that $\theta_1 = id$, then conditions
(1) and (2) (which read as condition (1) and
  (\ref{relation})) provide with the relations:
\begin{enumerate}
\item $-2i  \frac{\partial a}{\partial
  \bar{x}}=\omega(x,a(x))$
\item $\frac{\partial \omega}{\partial y}(x,a(x)) \, b(x) +2 i
  \frac{\partial b}{\partial \bar{x}}=0$.
\end{enumerate}
\par
Thus, the function $a$ is solution of a
$\bar{\partial}$-differential equation, and once a solution $a$
has been fixed, $b$ is also solution of a
$\bar{\partial}$-differential equation. This determines the
function $\theta$ and then the perturbation term in the expression
of $\tilde{\omega}$ (compared with $\omega$) which is
$\frac{\partial
  \theta_2}{\partial \bar{x}} =   \frac{\partial a}{\partial
  \bar{x}} + \frac{\partial b}{\partial \bar{x}} \, y$.
Thus all the possible functions $\theta$ corresponding to changes
of $F$, and so all the possible changes for $\omega$ are
determined.

\section{Invariants of the foliation}
The tensor $\Omega$ and the function $\omega$ measure the lack of
holomorphy of the function $F$. However, they depend on $F$. Since
we have noticed that one $F$ is \hle if and only if all $F \in
\mathcal{C}$ taking values in the image of $F$ are, this raises
the issue of the existence of an invariant of the foliation
measuring the lack of holomorphy of all the maps $F \in \mcl{C}$.
Looking at the relation (\ref{unepetiteeq}) expressing how varies
$\omega$ under change of $F$, we notice that by differentiating it
two times with respect to $y$, the variation in $\frac{\pp^2
\omega}{\pp y^2}$ doesn't contain any additive term anymore: \be
\label{secondderiveq} \frac{\partial^2
  \td{\omega}}{\partial y^2} = \frac{\partial^2
  \omega}{\partial y^2} \circ \theta \, \frac{\partial
  \theta_2}{\partial y} \, \overline{\frac{\partial
  \theta_1}{\partial x}}.
  \ee

\par

Thus, by considering the local differential forms on $X$ got by
pushing-forward the function $\omega$ (defined on $\D\times \C$),
and studying how they depend on $F$, we can check that the tensor
on $X$ got by differentiating them twice in an appropriate way
does not depend on $F$. Thus we are able to construct an invariant
of the foliation: a tensor, both global on $X$ and independent on
$F$, vanishing if and only if the maps $F$ are \hl, as wanted.

\subsection{Constructions of differential forms/tensors}
\label{paragdef} The map $F$ being a local diffeomorphism, there
exists, on a neighborhood $V_F$ of the point $z_0:=F(0,0)$, a
local inverse for $F$, taking values in a neighborhood $U_0$ of
the point $(0,0) \in \D \times \C$. It will be (abusively) denoted
by $F^{-1}: V_F \ra U_0$, and its coordinate-functions will be
denoted by $F^{-1}=((F^{-1})_1, (F^{-1})_2)$.\\

Thus the map $F$ provides us with a local trivialization of the
bundle $T\mathcal{F}$ on $V_F$:
\begin{equation*}
\phi_F~:
\begin{cases} T\mathcal{F}_{|V_F} & \lgra  V_F \times \C \\
  \xi_z & \lgra  \left(z, \pi_2(\dd F_z^{-1}(\xi_z)\right)=\left(z, \dd
  (F^{-1})_2 (\xi_z)\right)
\end{cases}
\end{equation*}
(with $\pi_2$ the second coordinate-function). Possibly
restraining the transverse disk $D$, let us assume the
neighborhood $V_F$ contains $F(\D\times \{ 0 \})$. Moreover, in
the following, we can assume that there exists a holomorphic
foliated chart on
$V_F$ (possibly restraining $V_F$) $\psi: V_F \ra \D \times \D.$\\

This being done, the tensor $\Omega$ can be read on $V_F$ as a
tensor taking values in ${T\mathcal{F}}_{|V_F}$ (from now on the
restriction to $V_F$ will always be implicit) :
\begin{equation*}
\Omega_0= (F_*)\circ \Omega \circ (F^{-1}_*): \begin{cases} TX & \lgra T\mathcal{F} \\
Z & \lgra JZ-F_*(J_0F^{-1}_* Z). \end{cases}
\end{equation*}
\noindent Then composing by $\phi_F$ we get a $1$-form $\lambda_F$
on $V_F$:
\begin{equation*}
\lambda_F : \begin{cases} TX & \lgra \C \\
Z & \lgra \pi_2(F^{-1}_*(JZ)-J_0(F^{-1}_* Z))=
\pi_2(\Omega(F^{-1}_*
 Z)) =(F^{-1})^* \Omega_2 \, (Z) \end{cases}
\end{equation*}
with $\Omega_2= \pi_2 \circ \Omega$ a $1$-form on $\D \times \C$,
which reads by definition of $\omega$ as $\Omega_2 \begin{pmatrix}
v_1 \\ v_2 \end{pmatrix}=\omega \, \bar{v_1}$. This form
$\lambda_F$ is zero along the leaves and $\lambda_F(JZ)=-J_0
\lambda_F(Z)$, that is to say $\lambda_F$ is a $(0,1)$-form. By
construction, $\lambda_F$ is zero if and only if $\Omega$
vanishes, thus if and only if $F$ is holomorphic on $V_F$.
Furthermore, we can point out that $\lambda_F = - 2 i\,
\bar{\partial}((F^{-1})_2)$ with the aforementioned notation.\\

Finally, by definition of $\Omega$, for every vector field $Z$,
$\lambda_F(Z)= \pi_2(\Omega(F^{-1}_* Z))=\omega(F^{-1}(z))\
\overline{\dd (F^{-1})_1 (Z)}= \omega(F^{-1}(z))\,\epsilon(Z)$,
with $\epsilon$ the $1$-form on $V_F$, $\epsilon(Z)=\overline{\dd
(F^{-1})_1 (Z)}.$ Let us notice that the form $\epsilon$ is equal
to $(F^{-1})^* \delta$, with $\delta$ the $1$-form on $F^{-1}(V_F)
\subset \D \times \C$ defined by $\delta(V)= \overline{\pi_1(V)}$
(with $\pi_1$ the projection on the first coordinate of $\C \times
\C$), \emph{i.e.} $\delta= \dd \bar{x}$ (if $(x,y)$ is the
standard coordinates system on $\D \times \C$). Thus $\dd
\delta=0$ and $\dd \epsilon=0$. So, the $1$-form $\epsilon$ on
$V_F$ is a basic $1$-form of the foliation: for any $Y \in T
\mcl{F}$, $\epsilon(Y)=0$ and $i_Y \dd \epsilon=0$.\\

\par

Let us now consider the $2$-form $\dd \lambda_F$ on $V_F \subset
X$, $\dd \lambda_F~: TX \times TX  \ra \C$. Since $\dd=\partial
+\bar{\partial}$ and $\bar{\partial}^2=0$, it satisfies: $\dd
\lambda_F= 2 \partial \bar{\partial} (F^{-1})_2$ on $V_F$ and $\dd
\lambda_F \in \Lambda^{1,1} X$. In regards of the expression of
$\lambda_F$, $\dd \lambda_F (Z,Y)= \dd \omega (\dd F^{-1}(Z)) \,
\epsilon(Y)-  \dd \omega(\dd F^{-1}(Y)) \, \epsilon(Z)) + \omega
\circ F^{-1} \, \dd \epsilon (Z,Y)$. Thus considering $Y\in T
\mcl{F}$, $\dd \lambda_F(Z)= - \frac{\partial \omega}{\partial
y}\, (F^{-1}(z)) \dd (F^{-1})_2(Y) \,\epsilon (Z)$.\\

Thus, restricted to $TX \times T\mcl{F}$, $\dd \lambda_F$ defines
a tensor:
\begin{equation*}
A_F~: \begin{cases} TX \times T\mathcal{F} & \lgra T\mathcal{F}  \\
(Z,Y) & \lgra {\phi_F}^{-1}(\dd\lambda_F(Z,Y)). \end{cases}
\end{equation*}
From the properties of $\lambda_F$, we deduce that $A_F$ is zero
on $T\mathcal{F} \times T\mathcal{F}$. Moreover $A_F$ can be
written in terms of $\omega$ : $A_F (Z,Y) = - \frac{\partial
\omega}{\partial y} (F^{-1}(Z)) \, \epsilon(Z) \, Y = -
\frac{\partial \omega}{\partial y} (F^{-1}(Z)) \,\overline{\dd
(F^{-1})_1 (Z)}\, Y$.\\

 \par It is is
$\C$-antilinear with respect to the first variable and $\C$-linear
with respect to the second one. Furthermore, $A_F$ is zero if and
only if $\frac{\partial \omega}{\partial y}$ is zero on $V_F$, or
equivalently if and only if $\omega=0$ on $ \D \times \C$
\emph{i.e.} if and only if $F$ is \hle on $\D \times \C$. In view
of the above expression of $\dd\lambda_F$, one can note that $F$
is \hle if and only if  $\partial \bar{\partial}(F^{-1})_2=0$ on
$TX
\times T \mathcal{F}$ {\it (harmonic condition)}.\\

Finally, let us point out that $A_F$ can also be seen as a
$1$-form taking values in $\mathcal{E}nd (T\mathcal{F})$~($\simeq
\C$ since $dim_\C T \mathcal{F}=1$): \beq*
A_F: \begin{cases} TX \ra \mathcal{E}nd (T\mathcal{F}) \simeq \C\\
\ X \ra (Y \ra A_F(Z,Y))
\end{cases}
\eeq* This is a $1$-form of type $(0,1)$ and can be expressed in
function of $\omega$: $A_F(Z)= - \frac{\partial \omega}{\partial
y} (F^{-1}(z)) \ \overline{\dd (F^{-1})_1 (Z)}$, \emph{i.e.}
according to the previous notation :

\be \label{afexpress} A_F=- \frac{\partial \omega}{\partial y}
(F^{-1}(z)) \ \epsilon. \ee

\subsection{Modifications when one changes $F$ into $\tilde{F}$}
\label{paragchang}

As previously, we assume $F(\D \times \C)=\td{F}(\D \times
  \C)$ and we introduce the function $\theta$, given
  by lemma \ref{deftheta}, such that $F= \theta \circ
\td{F}$. Then, we prove:
\begin{prop}
If we consider $A_F$ and $A_{\tilde{F}}$ as $1$-forms, then on
$V_F \cap V_{\td{F}}$
\begin{equation*}
\frac{\partial \theta_2}{\partial y}(\tilde{F}^{-1}(z))\
A_{\tilde{F}}(Z) = \frac{\partial \theta_2}{\partial
y}(\tilde{F}^{-1}(z)) \  A_F(Z) - 2 i \frac{{\partial}^2
  \theta_2}{\partial y \partial \bar{x}}(\tilde{F}^{-1}(z))\, \overline{\dd
  ({\tilde{F}}^{-1})_1(Z)},
\end{equation*}
\end{prop}

\begin{proof} By definition,
\begin{equation*}
\begin{pmatrix} 0 \\ \lambda_{\tilde{F}}(Z) \end{pmatrix} =
  (\tilde{F}^{-1})_* (JZ)-J_0 (\tilde{F}^{-1})_* (Z)
\end{equation*}
\noindent So

\begin{equation*}
\theta_* \begin{pmatrix} 0 \\ \lambda_{\tilde{F}}(Z) \end{pmatrix}
=
  \begin{pmatrix} 0 \\ \lambda_F(Z) \end{pmatrix} + J_0
  \theta_*(\tilde{F}^{-1}_* Z)-\theta_* (J_0 \tilde{F}^{-1}_* Z)
\end{equation*}

\noindent That is, $\frac{\partial \theta_2}{\partial y}
(\tilde{F}^{-1}(z)) \, \lambda_{\tilde{F}}(Z)=\lambda_{F}(Z)+ 2 i
\frac{\partial
  \theta_2}{\partial \bar{x}}(\tilde{F}^{-1}(z))\,
\overline{\dd({\tilde{F}}^{-1})_1 (Z)}$.\\
\noindent Thus, differentiating and evaluating on  $X \in TX$ and
$Y \in
T\mathcal{F}$, we get \\
$\frac{\partial \theta_2}{\partial y}(\tilde{F}^{-1}(z)) \, \dd
\lambda_{\tilde{F}}(Z,Y)=\dd \lambda_{F}(Z,Y)- 2 i
\frac{{\partial}^2
  \theta_2}{\partial y \partial \bar{x}}(\tilde{F}^{-1}(z))\, \dd
({\tilde{F}}^{-1})_2(Y) \, \overline{\dd ({\tilde{F}}^{-1})_1(Z)}$
(since $\dd({\tilde{F}}_1)(Y)=0$). Finally, as $\phi_F \circ
\phi_{\tilde{F}}^{-1} (z,v)=(z, \frac{\partial \theta_2}{\partial
y}(\tilde{F}^{-1}(z))\, v)$, this implies:

\begin{equation*}
\label{changF1} \frac{\partial \theta_2}{\partial
y}(\tilde{F}^{-1}(z)) \  A_{\tilde{F}}(Z,Y) = \frac{\partial
\theta_2}{\partial y}(\tilde{F}^{-1}(z)) \  A_F(Z,Y) - 2 i
\frac{{\partial}^2
  \theta_2}{\partial y \partial \bar{x}}(\tilde{F}^{-1}(z)) \, \overline{\dd
  ({\tilde{F}}^{-1})_1(Z)} Y.
\end{equation*}

\noindent Therefore, if we see $A_F$ as a $1$-form taking values
in $ \mathcal{E}nd (T\mathcal{F})$: $A_F: TX \ra  \mathcal{E}nd
(T\mathcal{F}) \simeq \C$, this relation reads as:
 \begin{equation}
\label{changF} \frac{\partial \theta_2}{\partial
y}(\tilde{F}^{-1}(z)) \  A_{\tilde{F}}(Z) = \frac{\partial
\theta_2}{\partial y}(\tilde{F}^{-1}(z)) \  A_F(Z) - 2 i
\frac{{\partial}^2
  \theta_2}{\partial y \partial \bar{x}}(\tilde{F}^{-1}(z))\, \overline{d
  ({\tilde{F}}^{-1})_1(Z)}.
\end{equation}
\end{proof}

Let us note that in regard with (\ref{relation}), this can be
written: \ba*
 A_{\tilde{F}}(Z) & =
 A_F(Z) - \frac{\partial
  \omega}{\partial y} (F^{-1}(\tilde{\pi}(z))) \,
 \overline{d({\tilde{F}}^{-1})_1(Z)}\,  \overline{\frac{\partial
 \theta_1}{\partial x}}(\tilde{F}^{-1}(z))\\
 & = A_F(Z) - \frac{\partial
  \omega}{\partial y} (F^{-1}(\tilde{\pi}(z))) \,
 \overline{d({F}^{-1})_1(Z)},
\ea* \noindent where $\tilde{\pi}$ is the projection of $V_F$ on
$\tilde{F}( \D \times \{ 0 \})$ along the leaves (see Appendix):

\begin{equation*}
\tilde{\pi}~: \begin{cases} V_F & \lgra V_F \\
z=\tilde{F}(x,y) & \lgra \tilde{F}(x,0). \end{cases}
\end{equation*}

Thus, in sight of the expression of $A_F$ in term of $\omega$, it
follows that:

\begin{equation*}
 A_{\tilde{F}}(Z) =
 A_F(Z) - {\tilde{\pi}}^* A_F (Z)
\end{equation*}

So, the perturbation term in the expression of $A_F$ under change
of the map $F$ depends only on the values of $A_F$ along the
holomorphic disk $\td{D}$. As a consequence, if $A_F$ identically
vanishes, then so does $A_{\tilde{F}}$; and we are able to recover
the statement that if $F$ is \hle then $\td{F}$ is too (on
$\td{F}^{-1} (F(\D \times \C))$.\\
\par
These differential forms associated to $F$ have been constructed
quite naturally. We would like to better understand their meaning.
In the next section we shall explain how this invariant $A_F$ is
closely related with the $\bar{\pp}$-connection on $T\mathcal{F}$.

\subsection{Connections on $T\mathcal{F}$}

The canonical connection $\bar{\partial}_0$ of the trivial bundle
with fiber $\C$ on $U_F$ is defined by $\bar{\partial}_0 (Z,
\xi)=\bar{\partial} \xi (Z)$. From this connection, the map $F$
determines a $\bar{\pp}$-connection on the bundle $T
\mathcal{F}_{|V_F}$: $\bar{\partial}_F=(\phi_F)^{-1}\circ
\bar{\partial}_0 \circ \phi_F$, which means for any section $\xi
\in T\mcl{F}$ and any vector field $Z \in TX$,
$(\bar{\partial}_F)_Z
(\xi)=(\phi_F)^{-1} (\bar{\partial}( \phi_F (\xi))(Z))$.\\

\par Then $\nabla_F =\frac{i}{2}
A_F+\bar{\partial}_F$ -- \emph{i.e.} for any section $\xi \in
T\mcl{F}$ and any vector field $Z \in TX$, $(\nabla_F)_Z(\xi)
=\frac{i}{2} A_F (Z) \xi+ (\bar{\partial}_F)_Z(\xi)$ -- defines a
$\bar{\pp}$-connection on the bundle ${T\mathcal{F}}_{|V_F}$.

This is a $\bar{\pp}$-connection and:
\begin{prop}
The connection $\nabla_F =\frac{i}{2}
  A_F+\bar{\partial}_F$ is independent from the choice of $F$. More precisely,
this is the canonical $\bar{\pp}$-connection on the \hle bundle
$T\mathcal{F}$ -- that is to say the one defining the \hle
structure of $T\mathcal{F}$.
\end{prop}
\begin{proof}
Using that $(\phi_F) \circ \phi_{\tilde{F}}^{-1} (z,v)=(z,
\frac{\partial \theta_2}{\partial y} v)$, we check that

\beq* (\bar{\partial}_F)_Z(\xi)= \phi_F^{-1} \left( \bar{\pp}
\left( \f{\pp \theta_2}{\pp y} \circ \td{F}^{-1}  \
\phi_{\td{F}}(\xi) \right)(Z) \right) \eeq* and then since
$\phi_F^{-1}(v) = \phi_{\td{F}}^{-1} \left(\f{1}{\f{\pp
\theta_2}{\pp y}} \, v\right)$,

\beq* (\bar{\partial}_F)_Z(\xi)= (\bar{\partial}_{\td{F}})_Z(\xi)
+ \f{1}{\f{\pp \theta_2}{\pp y}} \f{\pp^2 \theta_2}{\pp x \, \pp
y} \circ \td{F}^{-1} \, \overline{\dd (\td{F}^{-1})_1(Z)} \ \xi
\eeq* \noindent This inequality and (\ref{changF}) prove that
$\nabla_F=\nabla_{\tilde{F}}$.\\

\noindent Moreover, if we introduce a foliated chart $\psi~: V_F
\ra \D \times \D$, the canonical $\bar{\pp}$-structure of the \hle
bundle $T\mcl{F}$ is $(\bar{\pp}_\mcl{F})_Z (\xi) = \psi_*^{-1}
\left(\bar{\pp}\left( \psi_*(\xi) \right)(Z)\right)$. \\
The same kind of calculations as in section \ref{paragchang} --
but using $\beta = F^{-1} \circ \psi$ instead of $\theta $ --
allows to read all the forms defined in \ref{paragdef} in this
chart and we get:
\begin{equation*}
\frac{i}{2} d\lambda_F(Z,Y)= - \frac{\partial^2}{\partial y
  \partial \bar{x}}(\psi(z))\, \overline{d\psi_1(Z)} d\psi_2(Y).
\end{equation*}
Then, using the same argument as the one used above to prove that
$\nabla_F= \nabla_{\td{F}}$ -- but again using $\beta$ instead of
$\theta$ -- leads to $\nabla_F= \bar{\pp}_\mcl{F}$, the canonical
$\bar{\pp}$-connection on the \hle bundle $T\mathcal{F}$.
\end{proof}

\subsection{Invariant tensor field $\Gamma$}
\label{invt} In light of our previous study, we can construct a
tensor field, independent of $F$, an invariant of the foliation
that vanishes if and only if the maps $F \in \mathcal{C}$ are \hl,
\emph{i.e.} if and only if there exist some \hle maps $F:\D\times
\C \ra X$ compatible with the foliation. This answers our initial
question of the existence of holomorphic cylinders along the
foliation.

\begin{defi-thm}
Let  $\Gamma_F=(\dd A_F)_{| TX \times T\mathcal{F}}$~:
\begin{equation*}
{\Gamma_F}_{|x}:\begin{cases} (T_x X \times T_x \mathcal{F}) &
\lgra
   \C \\
(Z,Y) & \lgra \dd A_F(Z,Y) \end{cases}, \ \forall x \in V_F
\end{equation*}
It does not depend on $F$. Thus this defines a tensor $\Gamma$:
$\Gamma=\Gamma_F$ on each $V_F$. This is a tensor on $X$ which is
associated with the foliation, $\C$-antilinear with respect to the
first variable and $\C$-linear with respect to the second one. It
vanishes if and only if the maps $F$ are holomorphic. More
precisely, for any $F \in \mathcal{C}$, $\Gamma$ vanishes on
$F(\D\times \C)$
if and only if $F$ is \hl.\\
Moreover, $\Gamma$  is a leaf-wise \hle section of the line bundle
$E=\Lambda^{(0,1)} T^* \mathcal{N} \otimes \Lambda^{(1,0)} T^*
\mathcal{F}$, \hle above $(X, \mathcal{F})$.
\end{defi-thm}
\begin{proof}
As seen previously, $A_F$ is a $\C$-valued $1$-form on $V_F$. So,
$\dd A_F$ is a $2$-form. Since $\f{\pp \theta_2}{\pp y}$ is
independent of $y$ and $\epsilon= \dd (\td{F}^{-1})$ is a basic
$1$-form, differentiating (\ref{changF}) proves that, if
restricted to $TX \times T\mathcal{F}$, $\dd A_F$ does not depend
on $F$. Thus, we can define $\Gamma$ on $X$.\\
More precisely, according to the expression (\ref{afexpress}) for
$A_F$ (and since $\epsilon$ is a basic form), $\Gamma$ can be
written locally as $\Gamma(Z,Y)= \dd \left(\f{\pp^2 \omega}{\pp
y^2}\right) (\dd (F^{-1}(Y)) \ \epsilon(Z)$, and so:
\begin{equation}
\label{espressgamma} \Gamma(Z,Y)=  \f{\pp^2 \omega}{\pp y^2}
(F^{-1}(z)) d (F^{-1})_2(Y) \, \overline{d(F^{-1})_1(Z)}.
\end{equation}

\noindent Moreover, $\Gamma$ is antilinear with respect to the
first variable, linear with respect to the second one. Thus
$\Gamma$ belongs to $E=\Lambda^{(0,1)} T^* {X} \otimes
\Lambda^{(1,0)} T^* \mathcal{F}$. Furthermore, from the results on
$A_F$, we know that if the maps $F$ are \hl, then $\Gamma$
vanishes on $F(\D \times \C)$.

\noindent Conversely, if $\Gamma$ vanishes on $F(\D \times \C)$,
it follows that $\frac{\partial^2  \omega}{\partial y^2} =0$ (on
$F^{-1}(V_F)$ and so on $\D \times \C$), and consequently,
$\omega$ being \hle along the leaves, $\frac{\partial
\omega}{\partial y}$ is constant on each leaf. However, as this
function vanishes on $\D \times\{ 0\}$, it is identically zero.
Thus, $\omega$ vanishes identically.
Finally the maps $F$ are \hle if and only if $\Gamma=0$.\\

\noindent Finally, since $\Gamma$ vanishes on $T\mathcal{F}\times
T\mathcal{F}$, it can be considered as a form on
$T\mathcal{N}\times T\mathcal{F}$, which is $J$-antilinear with
respect to the first variable and $J$-linear with respect to the
second one. So $\Gamma$ is a section of the leaf-wise \hle bundle
$E=\Lambda^{(0,1)} T^* \mathcal{N} \otimes \Lambda^{(1,0)} T^*
\mathcal{F}$ above $(X, \mathcal{F})$. Furthermore, in sight of
(\ref{espressgamma}), it reads in the \hle foliated chart $\psi:
V_F \ra \D \times \D$ : $\left((\psi^{-1})^* \Gamma
\right)_{(x,y)} (v_x, v_y) = \f{\pp^2 \omega}{\pp y^2} (\beta (z))
\, \f{\pp \beta_2}{\pp y}(x,y) \, \f{\pp \beta_1}{\pp x}(x) \, v_y
\, \overline{v_x},$ with $\beta =(\beta_1, \beta_2)= F^{-1} \circ
\psi^{-1} : \D \times \D \ra \D \times \C$.\\
Thus, we notice that $\Gamma$ is \hle along the leaves. \\
The same reasoning also proves that the bundle $E$ is \hle along
the leaves.
\end{proof}

Therefore the issue of the existence of some \hle map $F:\D \times
\C \ra X$ that is compatible with the foliation can come down to
the study of the existence of leaf-wise
\hle sections of the bundle $E$ above the foliation $\mathcal{F}$.\\

Let us make a small remark: we could have deduced that $\Gamma_F$
is independent of $F$ from its expression (\ref{espressgamma}) in
function of $F$, using the equality (\ref{secondderiveq}).\\

Another remark can be made: at the end of the proof above, we have
expressed $\Gamma$ in the foliated chart $\psi$ in terms of
$\omega$ and $\beta$. In fact, it can be expressed only in terms
of $\beta$. Indeed, writing the holomorphy of $\beta$ between $(\D
\times \D, J_0)$ and $(\D \times \C, J_0 + \Omega)$ leads us to: $
- 2 i \, \f{\pp \beta_2}{\pp \bar{x}}= \omega \circ \beta \,
\overline{\f{\pp \beta_1}{\pp x}}.$ So if we suppose, without any
loss of generality that $\beta_1=id$, then \beq* \omega \circ
\beta = - 2 i \, \f{\pp \beta_2}{\pp \bar{x}}. \eeq* (This could
also have been deduced from the relation (\ref{omegalapharel}.)
Differentiating twice this equation, we get: $\f{\pp^2 \omega}{\pp
y^2} \circ \beta \, \left(\f{\pp \beta_2}{\pp y} \right)^2=
-\left( \f{\pp \omega}{\pp y} \circ \beta \, +\, 2i \f{\pp^3
\beta_2}{\pp \bar{x} \pp y^2} \right).$ Thus if $\psi$ is a
foliated chart sending $F(\D \times \{ 0\}$ to $\D \times \{0\}$
(identically, \emph{i.e} via $\alpha_1=id$), then on this disk,
$\Gamma$ reads in this chart:

\beq* \left( (\psi^{-1})^* \Gamma \right)_{(x,0)} (v_x, v_y)= - 2
i \,  \left( \f{\pp \beta_2}{\pp y}(x,0) \right)^{-1} \, \f{\pp^3
\beta_2}{\pp \bar{x} \pp y^2}(x,0). \eeq*

\subsection{Generalization to higher dimensions}
\label{higherdim} Our analysis, written in the case of a complex
dimension $2$ manifold $X$, can be straightforwardly generalized
to the case of higher dimensions. Let $X$ be a manifold of complex
dimension $n$ provided with a $1$-dimensional foliation $\mcl{F}$
with parabolic leaves. We can consider the cylinder-maps $F$ from
$\D^{n-1} \times \C$ into $X$, constructed along a \hle transverse
disk $D \simeq \D^{n-1}$ and in the direction of a trivialization
$t$ of the line bundle $T\mcl{F}$ along $D$. Then, as previously,
we can consider the tensor $\Omega$ on $\D^{n-1} \times \C$
defined by $J_0+ \Omega = F^* J$. It can be written:
\begin{equation*}
\Omega_{(x,y)} \begin{pmatrix} v_x \\ v_y
\end{pmatrix}=\begin{pmatrix} 0 \\ \omega_{(x,y)} (v_x) \end{pmatrix}.
\end{equation*}
with for any $(x,y) \in \D^{n-1} \times \C$, $\omega_{(x,y)}$ a
$\C$-antilinear $1$-form on $\C^{n-1}$.\\
\par
Considering two functions $F$ and $\td{F}$, we can, as before,
construct a map $\theta~: (x,y) \in \D^{n-1} \times \C \ra
(\theta_1(x), \theta_2(x,y)) \in \D^{n-1} \times \C$, with
$\theta_1$ a bi-holomorphism of $\D^{n-1}$, and for any $x$,
$\theta_2(x,.)$ a bi-holomorphism of $\C$, which satisfies
$\td{F}=F \circ \theta$. Then, as previously we check, \beq*
{\frac{\partial
  \theta_2}{\partial y}}\tilde{\omega}(x,y)=
 \left( \frac{\partial  \theta_1}{\partial x}\right)^* \
 \omega(\theta(x,y)) + 2 i \frac{\partial \theta_2}{\partial
 \bar{x}},
\eeq* \noindent where $ \left( \frac{\partial  \theta_1}{\partial
x}\right)$ is a linear endomorphism of $\C^{n-1}$ and where $*$
denotes its
conjugate (for the canonical hermitian product on $\C^{n-1}$).\\
\par Then, in the same way we can define, from each $F$, a local differential
form $A_F: TX \ra End(T \mcl{F}) \simeq \C$, and an invariant
$\Gamma$:

\begin{defi-thm}
Let  $\Gamma_F=(\dd A_F)_{| TX \times T\mathcal{F}}$~:
\begin{equation*}
{\Gamma_F}_{|x}:\begin{cases} T_x X \times T_x \mathcal{F} & \lgra
  \mathcal{E}nd(T\mathcal{F}) \simeq \C \\
(X,Y) & \lgra \dd A_F(X,Y) \end{cases}, \ \forall x \in V_F
\end{equation*}
It does not depend on $F$. Thus this defines a tensor $\Gamma$:
$\Gamma=\Gamma_F$ on each $U_F$. This is a tensor on $X$ which is
associated to the foliation, $\C$-antilinear with respect to the
first variable and $\C$-linear with respect to the second one. It
vanishes if and only if the maps $F$ are holomorphic. More
precisely, for any $F \in \mathcal{C}$, $\Gamma$ vanishes on
$F(\D\times \C)$
if and only if $F$ is \hl.\\
 Moreover, $\Gamma$  can be seen as a leaf-wise
\hle section of the vector bundle $E=\Lambda^{(0,1)} T^*
\mathcal{N} \otimes \Lambda^{(1,0)} T^* \mathcal{F}$ above $(X,
\mathcal{F})$.
\end{defi-thm}

\par \vspace*{.5cm} After this general study, we would like to analyze a few
particular cases in which some more specific results can be
proven.
\section{Particular cases and examples}

Let us recall that in the case of Stein manifolds, it has been
proven that all the maps $F \in \mathcal{C}$ are \hle (that is
$\Gamma=0$), \cite{nishino} (see also \cite{il1}, \cite{il2}). Our
study allows us to show that, in the case of a foliation with
compact leaves too, all the maps $F$ are \hle (\emph{i.e.}
$\Gamma=0$).

\subsection{Foliation with compact leaves, and holonomy groups}
\label{compactleaves}

In the particular case of foliations with compact leaves, several
known results of stability are at our disposal.\\
First, let us recall the equivalence (see \cite{god}):

\begin{prop}
Let $\mcl{F}$ be a foliation on a manifold $X$ whose leaves are
all compact. Then the following properties are all equivalent:
\begin{enumerate}
\item each leaf has a finite holonomy group

\item each leaf has a fundamental system of saturated open
neighborhoods

\item the space of leaves $M / \mcl{F}$ is Hausdorff

\item the saturation of a compact set is a compact set

\end{enumerate}
Such a foliation is called {\it stable}.
\end{prop}

It has been proven that these properties are satisfied in some
special cases. Notably in the case of holomorphic foliation (see
\cite{holmann} and \cite{holmann2}, compared also with
\cite{edsull}):
\begin{prop}
Let $X$ be a complex manifold provided with a \hle foliation with
compact leaves. If $X$ is K\"ahler, or if the foliation has
complex codimension $1$, then the foliation is stable.
\end{prop}

Thanks to these results, we can deduce:

\begin{prop}
\label{compleavesprop} Let $X$ be a complex manifold provided with
a \hle foliation of dimension $1$, such that $X$ is either
K\"ahler, or has complex dimension $2$.\\
If $D$ is a transverse \hle disk whose saturation is foliated by
compact leaves, then all the maps $F \in \mathcal{C}$ constructed
along
this disk are \hl.\\
In other words, the invariant $\Gamma$ associated to the foliation
vanishes on the saturation of $D$.\\
As a consequence, if $X$ is foliated by compact leaves then
$\Gamma$ identically vanishes on $X$.
\end{prop}

\begin{proof}
Let us denote the saturation of $D$ by $U_0$.  And let us consider
a map $F \in \mathcal{C}$, constructed along the disk $D$, with
$t$ being the fixed trivialization of $T\mcl{F}_{| D}$. Possibly
changing the projection $\psi_x~: \td{L_x}\simeq \C \ra L_x$, we
can assume that for every $x \in \D$, $t_x=1$ in $\td{L_x}$. Then,
for all $(x,y) \in \D \times \C$, (keeping the notation from
section \ref{constructF} $F_x=\psi_x \circ F^0_x$) $F^0_x(y)=y \in
\td{L_x}$. Moreover, $L_x$ is written as ${\C}{/ G_x}$, with $G_x=
\gamma_1(x) \Z \oplus \gamma_2(x) \Z$ being its
fundamental group.\\

\noindent Let us fix a distinguished open set $U$ containing $D$
and a foliated chart $\psi~: U \ra \D \times \D$ such that $\psi
(D)=\D \times \{ 0 \}$. Then, we can consider the projection on
the transverse $D$ along the leaves (of $\mcl{F}_U$) $\pi_0 : U
\ra D$ (see Appendix \ref{rappel}).\\

\noindent Let us denote $x_0=F(0,0)$. For every $\gamma$ in the
fundamental group $G_0=G_{x_0}$, $F(0, \gamma)=x_0$. Thus, for
every $x \in D$ close enough to $x_0$, $F(x, \gamma) \in U$. So
there exists a unique $\alpha_\gamma(x) \in \D$ such that
$F(\alpha_\gamma(x),0)=\pi_0(F(x, \gamma))$. The map
$\alpha_\gamma$, thus defined, is the holonomy map associated with
$\gamma$. As explained in the Appendix, this map $\alpha_\gamma$
is holomorphic. First, let us show that
\begin{lemma}
There exists a disk $D_0 \subset D$ and $N \geq 1$ such that
$\alpha_\gamma$ is well-defined on $D_0$ and takes values into
$D_0$. Moreover, $\alpha_\gamma^N=Id$ on $D_0$.
\end{lemma}
\begin{proof}

According to the stability results in the case of foliations whose
all leaves are compact, the foliation induced on $U_0$ is stable.
So the leaf $L_{x_0}$ has a fundamental system of
saturated open neighborhoods.\\

\noindent Since $L_{x_0}$ is compact, one can consider a sub-disk
$D' \subset D$ such that $L_{x_0} \cap \bar{D'} = \{x_0\}$. Then
one can consider a saturated open set $V_0$ included in the
saturation of $D'$ (and so included into $U_0$). Since $L_{x_0}
\cap \bar{D'} = \{x_0\}$, $V_0 \cap D'$ is connected. Let us
denote this by $D_0$. By construction, for any $x \in D_0$,
$\alpha_\gamma(x) \in D_0$. So the map $\alpha_\gamma :D_0 \ra
D_0$ can be defined. Moreover, this is a bi-holomorphism of the
disk $D_0$
whose inverse is $\alpha_{-\gamma}$.\\

%In order to prove this lemma, we first notice that, since
%$L_{x_0}$ is compact, we can restrict $D$ to $D' \subset D$ so
%that $L_{x_0} \cap \bar{D'} = \{x_0\}$. Furthermore, $D'$ can be
%chosen small enough so that $\alpha_\gamma$ and $\alpha_{-\gamma}$
%could be defined on $D'$. Then, for any small neighborhood $V$ of
%$x_0$ in $D'$, we define the sets $\displaystyle U_0:=\cup_{x \in
%V} L_x$ and $V_0:=U_0 \cap D'$. Because of the compactness of the
%leaves, $\cup_{x \in \bar{V}} L_x$ is compact.\\
%Moreover, for $V$ chosen small enough, $V_0$ is strictly included
%in $D'$ and there even exists $\epsilon>0$ such that $B(V_0,
%\epsilon) \subset D'$. Indeed, otherwise, one would get sequences
%of points $x_n \ra x_0$ and $z_n \in L_{x_n} \cap D'$ such that
%$z_n \ra \pp D'$; then using the compactness of the leaves, we
%could extract a sub-sequence of $(z_n)$ converging to $z_0 \in
%L_{x_0} \cap \pp D'$, which is impossible.\\

%Then, let us note that this set $V_0$ is connected (since $L_{x_0}
%\cap \bar{D'} = \{x_0\}$). Thus, choosing $D_0=V_0$, the map
%$\alpha_\gamma :D_0 \ra D_0$ can be defined.

\noindent Then for any $x \in D_0, \   \{\alpha_\gamma^n (x), n
\in \N \}$ is included in $L_x \cap \bar{D_0}$. The leaf $L_x$
being compact, this is a compact, discrete set, and therefore
finite. So, there must exist $N \geq 1$ such that
$\alpha_\gamma^{N_x} (x)=x$. Therefore,
\begin{equation*}
\dsps D_0 \subset \cup_{n \geq 1} \{x \in \bar{D_0} \textrm{ such
that } \alpha_\gamma^{n} (x)=x \}.
\end{equation*}
Since each of the sets of this union is closed, necessarily, one
among them has a non-empty interior. So, there exists a $N$ such
that $\alpha_\gamma^{N} (x)=x$ on an open set $V \subset D_0$. As
$\alpha_\gamma^{N}$ is holomorphic, we can conclude that
$\alpha_\gamma^{N}=Id$ on $D_0$.
\end{proof}

\noindent This being proven, if the fundamental group $G (x_0)=
G_0=\gamma_1 \Z \oplus \gamma_2 \Z$, we can fix a $N$ such that
for $i=1,2$, $\alpha_{\gamma_i}^{N}=Id$ on $D_0$ (possibly
restraining it to a smaller disk). Then, $\pi_0(F(x, N
\gamma_i))=F(x,0)$ on $D_0$. So, the functions $N \gamma_i(x_0)$
can be continuously extended on $D_0$ by a map $\beta_i(x) \in
G_x$. Let $G'_x$ be the subgroup of $G_x$ generated by these maps
$\beta_i(x)$.Then, the functions $F(x,.)$ are
$G'_x$-periodical for $x \in D_0$.\\

\noindent Thus for any fixed family $\gamma(x) \in G'_x$,
continuous with respect to $x$, and for all $(x,y) \in \D \times
\C$, $F(x, y+\gamma(x))= F(x,y)$. Therefore, we check, either by
differentiating the equality $F(x, y+\gamma(x))=F(x,y)$ and
applying $J$, or by applying the results of section
\ref{paragchang} to the map $\td{F}(x,y)=F(x,y+\gamma(x))$ (that
belongs to $\mathcal{C}$ and for which $\theta (x,y)= (x,
y+\gamma(x))$):
\begin{equation*}
\omega(x,y+\gamma(x))=\omega(x,y)  -2i\,\f{\pp
  \gamma}{\pp \bar{x}}
\end{equation*}

\noindent This implies $\f{\pp \omega}{\pp y}(x,\gamma(x))= \f{\pp
\omega}{\pp y}(x,y)$. Thus for any fixed $x \in D_0$, the function
$\f{\pp \omega}{\pp y}(x,.)$ is \hle on  $\C$ and
$G'_x$-periodical. Therefore it is constant. As it vanishes on 0,
it vanishes identically on $\C$. Finally $\omega(x,.)$ vanishes on
$\C$ for any $x \in D_0$, neighborhood of $x_0$. This can be
applied around any other point in $D$. Thus, $\omega$ vanishes
identically on the saturation of $D$ and $F$ is \hl.
\end{proof}

Let us note that $\omega$  measures the lack of holomorphy of the
map $x \ra G'_x$. Indeed, if we fix any continuous map on $\D$
$\gamma(x)=j\beta_1(x)+k\beta_2(x) \in G'_x$, then $x \ra
F(x,\gamma(x))=F(x,0) \ \in L_x$ is $J$-\hl. Thus, $x \ra
(x,\gamma_x)$ is $J'$-holomorphic, which reads as:
 \beq*  -2i\, \f{\pp \gamma}{\pp
  \bar{x}}=\omega(x,\gamma(x)).
\eeq*

\par So, as a side product, through this proof, we get that the
groups $G'_x$ depend holomorphically on $x \in D$ \emph{i.e.} on
the transversal disk to the foliation. More precisely,

\begin{lemma}
Let $X$ be satisfying the same assumptions as in the previous
proposition. Let us suppose that along a transverse \hle disk $D$,
all the leaves are compact. Then if for any $x$, $G_x$ denotes the
fundamental group of the leaf $L_x$, any continuous family
$\gamma(x) \in G_x$, $x \in D$ depends holomorphically on $x$.
\end{lemma}

\par
One of the main points of these results is the finiteness of the
holonomy groups of the leaves. In the case of a foliation
$\mcl{F}$ possessing one compact leaf, two cases can be
distinguished depending on the holonomy group of the leaf. If the
holonomy group of the compact leaf is finite then we can use the
result (\cite{god} II.2.16.):
\begin{theorem}[Reeb]
A compact leaf from a foliation $\mcl{F}$ whose holonomy group is
finite, possesses a fundamental system of neighborhoods which are
saturated by compact leaves with finite holonomy groups.
\end{theorem}
In particular, if $L$ is a compact leaf with finite holonomy group
then there exists a small transverse disk $D$ around $L$ such that
all leaves crossing $D$ are compact. This result and the
proposition \ref{compleavesprop} allow us to deduce:

\begin{cor}
If $\mcl{F}$ is a \hle foliation of $X$ containing a compact leaf
$L_0$ with finite holonomy group, then there exists a \hle
cylinder $F: \D \times \C$ along the foliation with $F(0,0) \in
L_0$.
\end{cor}

If the holonomy group of the compact leaf $L_0$ is infinite then
the same reasoning does not apply. However in this case, we can
study how the holonomy group acts on $\omega$. Since $L_0$ is a
compact leaf whose universal covering is $\C$, $L_0$  is a torus
and $\pi_1(L_0) \simeq \Z^2$.\\

Let us consider $F \in \mcl{C}$, $F: \D \times \C \ra X$,
constructed along a transverse \hle disk $D_0$ and a fixed
trivialization $t$ of $ T \mcl{F}_{| D_0}$, and such that $x_0=
F(0,0) \in L_0$. Let $\omega$ be the function on $\D \times \C$
associated with $F$. Possibly by reparametrizing the universal
coverings $\td{L_x}$, we can assume that in $\td{L_x}$ $F(x,0)=0$,
\emph{i.e.} $F^0_x=0$ -- with the usual notations --, and $t_x =
1$, \emph{i.e} $(F^0_x)'(0)=1$. Then, these parametrizations of
$\td{L_x}$ being fixed, for any $x \in \D$, $F^0_x (y)=y \, \in
\td{L_x}$. And $L_0=L_{x_0} = \td{L_{x_0}} / G_0$ with $G_0=
\gamma_1 \Z \oplus \gamma_2 \Z= \pi_1(L_0)$ the fundamental group
of $L_0$.\\

In the following, to simplify the notations, for $x \in \D$, we
will denote the point $F(x,0) \in D$ by $x$ (\emph{i.e} we
identify $D$ with its \hle parametrization $x \ra F(x,0)$).\\

For any $\gamma_0 \in G_0$, $F( 0, \gamma_0)= F(0,0) = x_0$. One
can consider $\td{\gamma_0}$, the element of the holonomy group
associated with $\gamma_0$ (whose definition is explained Appendix
\ref{projhol}. Then (as in the proof of proposition
\ref{compleavesprop}) for any $x$ close to $x_0$ (identified with
$0$), there exists an element ${\gamma}(x)$ close to $\gamma_0$
such that $F(x, \gamma(x)) \in D$, thus defining a continuous map
$\gamma$ on a neighborhood $D_1$ of $0$, such that
$\gamma(0)=\gamma_0$. The map $\td{\gamma}$ defined by
$\td{\gamma}(x) = F(x, \gamma(x))$ is a local bi-holomorphism
around $0$ -- from $D_1$ to another neighborhood $D_2$ around $0$
--, whose derivative in $0$ is the holonomy element
$\td{\gamma_0}$.\\

Let us first note that since the map $x \ra F(x, \gamma(x))
=\td{\gamma}(x) \ (=F(\td{\gamma}(x),0))$ is \hl, the map $x \ra
(x, \gamma(x))$, defined on $D_1 \subset \D$ and taking values in
$\D \times \C$, is $(J_0, J')$-\hle (with $J'= J_0+ \Omega$). This
reads as

\be \label{toujourslameme} - 2 i \, \f{\pp \gamma}{\pp
\bar{x}}=\omega(x,\gamma(x))\ee

(The map $\omega$ measures the non-holomorphy of $\gamma$).\\

Let us consider the map $\td{F}$ in $\mcl{C}$ constructed along
$D_2 \subset D$, parametrized by $x \ra F(x, \gamma(x)) =
\td{\gamma}(x)$ ($=F(\td{\gamma}(x),0)$ with which it is
identified), $x \in D_1 \subset \D$, and with $\td{t}_x=
t_{\td{\gamma}(x)}$.\\

On one hand $\td{F}(x,0)=F(\td{\gamma}(x),0)$ and
$\td{t}_x=t_{\td{\gamma}(x)}$. So for any $x \in D_1$,
$\td{F}(x,y)= F(\td{\gamma}(x),y)$.\\

On the other hand, $F(x, \gamma(x))= \td{F}(x,0)$. So $\td{F}= F
\circ \theta$ with $\theta_1 = Id$ and $\theta_2(x,0) =\gamma(x)$.
Moreover, $\theta(x,.)$ is a bi-holomorphism of $\C$, and
consequently a degree-$1$-polynomial. So finally, $\theta(x,y)=(x,
\gamma(x)+ \phi(x) \, y)$. Moreover, since $\td{F}(0,0)=F(0,
\gamma(0))=F(0,0)$ and $\td{t}_0= t_0$, $\td{F}(0, .)= F(0,.)$ and
so $\phi(0)=1$.\\

If $\td{\omega}$ is the function associated with $\td{F}$,
 according to the relation (\ref{unepetiteeq}) in lemma
\ref{deftheta}, this leads to:

\be \label{facilemaisutile} \td{\omega}(x,y) =
\omega(\td{\gamma}(x),y)\, \overline{\f{\pp \td{\gamma}}{\pp
x}}(x). \ee

Moreover this also implies

\beq* \phi(x)\, \td{\omega}(x,y)=\omega(x, \gamma(x)+ \phi(x) \,
y) + 2i \left( \f{\pp \gamma}{\pp \bar{x}} + \f{\pp \phi}{\pp
\bar{x}}\, y \right),\eeq* which according to
(\ref{toujourslameme}) and (\ref{facilemaisutile}) leads to

\beq* \phi(x) \, \overline{\f{\pp \td{\gamma}}{\pp x}}(x) \,
\omega(\td{\gamma}(x),y) =\omega(x, \gamma(x)+ \phi(x) \, y) -
\omega (x, \gamma(x)) + 2 i \f{\pp \phi}{\pp \bar{x}}\, y. \eeq*

Differentiating this equality two times with respect to $y$, we
get:

\beq* \overline{\f{\pp \td{\gamma}}{\pp x}}(x)\, \f{\pp^2
\omega}{\pp y^2}(\td{\gamma}(x),y)= \f{\pp^2 \omega}{\pp
y^2}(x,\gamma(x)+ \phi(x) \, y) \, \phi(x).\eeq*

Thus at $x=0$,

\be \overline{\td{\gamma_0}} \, \f{\pp^2 \omega}{\pp y^2}(0,y) =
\f{\pp^2 \omega}{\pp y^2}(0,\gamma_0+ y),  \ee

\noindent with $\td{\gamma_0}$ the holonomy element associated
with $\gamma_0$. So to conclude, the holomorphic function
$h(y)=\f{\pp^2 \omega}{\pp y^2}(0,y)$ satisfies: for any $\gamma
\in G_0=\pi_1(L_0)= \gamma_1 \Z \oplus \gamma_2 \Z$, $h(y+
\gamma)= \overline{\td{\gamma}} \, h(y)$ with $\td{\gamma}$ the
holonomy element associated with $\gamma$. As a consequence, there
must exist some constants $C$ and $\rho$ such that $h(y)= C \,
\exp(\rho \, y)$ with $\exp(\rho \,
\gamma)=\overline{\td{\gamma}}$ for any $\gamma \in G_0$.\\

Thus, in the case where there exists a compact leaf $L_0$ whose
holonomy group $\td{G}_0$ is infinite, if there exists an
exponential map between the fundamental group $G_0 \subset \C$ of
the leaf and the holonomy group $\td{G}_0 \subset \C^*$,
$k:\begin{cases} G_0 & \ra \td{G}_0\\
\gamma & \ra  \overline{\exp(\rho \, \gamma)}
  \end{cases}$, then the function $\omega$ grows exponentially
  along $L_0$ : $\omega(y)= C \, \exp(\rho \, y)$.\\

Otherwise, $\omega$ vanishes identically on $L_0$.\\

\par If the manifold $M$ is compact -- even though the leaves are
not -- there is still a case where one can deduce that the
invariant $\Gamma$ vanishes, and therefore there is no \hle
cylinder along the foliation: the case where the bundle
$E=\Lambda^{(0,1)} T^* \mathcal{N} \otimes \Lambda^{(1,0)} T^*
\mathcal{F}$ above $(X, \mathcal{F})$ has a negative curvature.

\subsection{Case of bundle $E$ with negative curvature}

\begin{prop}
Let $X$ be compact. If the bundle $E=\Lambda^{(0,1)} T^*
\mathcal{N} \otimes \Lambda^{(1,0)} T^* \mathcal{F}$ above $(X,
\mathcal{F})$ has a strictly negative curvature along the leaves,
then the invariant $\Gamma$ vanishes identically and there is no
\hle cylinder along the foliation.
\end{prop}

\begin{proof}
Let us fix a Hermitian metric $g$ of $T\mcl{F}$ on $X$. This
provides each leaf with a Riemannian metric. If the curvature of
the bundle $E$ is strictly negative, then there exists a constant
$\epsilon>0$ such that $K < - \epsilon \ g.$\\

\noindent On one hand, the invariant $\Gamma$ is a leaf-wise \hle
section of the bundle $E$. So, if $\Gamma$ does not vanish
identically along a leaf, the curvature $K$ of $E$ along this leaf
is $K=\dd \dd^c \log |\Gamma|^2$ (except possibly on the numerable
set of points where $\Gamma=0$). Thus along such a leaf $K=\dd
\dd^c \log |\Gamma|^2 < - \epsilon \ g.$\\

\noindent On the other hand, the manifold $X$ being compact, the
curvature $K_g$ of $(T \mcl{F},g)$ is bounded by a constant on
$X$. In particular, there exists a constant $C$ such that along
the leaves
$ - C\, g <K_g< C \, g$.\\

\noindent If $\Gamma$ is not identically zero, let us consider a
leaf $L$ on which $\Gamma$ does not vanish identically. This leaf
(or its universal cover) is isomorphic to $\C$. Thus $\C$ is
provided with the metric $g= g_{|L}$, whose curvature satisfies
$K_g> -C \, g$. And we have a bundle $E_{|L}$ over $\C$, with a
\hle section $\Gamma$, whose curvature $K=\dd \dd^c \log
|\Gamma|^2 < - \epsilon \ g$. \\

\noindent Moreover, $K_g=\dd \dd^c \log |z|_g^2$ on $\C \exc \{
0\}$, with $z$ being a \hle section of $T\C$. So,

\beq* \dd \dd^c \log |\Gamma|^2 < \f{\epsilon}{C} \dd \dd^c \log
|z|_g^2,  \eeq* which reads as $\dd \dd^c \left(\log |\Gamma|^2 -
\f{\epsilon}{C} \log |z|_g^2 \right)<0.$ Thus $\phi=\log
|\Gamma|^2 - \f{\epsilon}{C} \log |z|_g^2$ is a sub-harmonic
function on $\C$ (with maybe a numerable set of points where
$\phi$ is not defined and in which $\phi$ goes to $- \infty$). In
particular, by the maximum principle, for any $r>0$, $\sup_{|z|=r}
\phi \geq \phi(0)$. Since $\Gamma$ does not vanish identically, we
can suppose that $|\Gamma(0)| \neq 0$. So, for any $r>0$, there
exists $z_r$ such that $|z_r|=r$ and $\phi(z_r) \geq \phi(0)=C_0$.
Finally, this means \beq* \log |\Gamma (z_r)|^2 \geq
\f{\epsilon}{C} \log r^2 + C_0. \eeq*

\noindent This implies that $|\Gamma (z_r)| \ra \infty$ when $r
\ra \infty$, which is not possible since $X$ is compact and so
$|\Gamma|$ is bounded on $X$. So $\Gamma$ has to vanish
identically.
\end{proof}

\par
In the general case, the invariant $\Gamma$ is not trivial: there
exist some foliated manifolds whose invariant $\Gamma$ does not
vanish and so, does not admit any holomorphic functions $F: \D
\times \C \ra X$ compatible with the foliation. Let us consider
the following example \footnotemark. \footnotetext{following
suggestions of Alexei Glutsyuk and Etienne Ghys}

\subsection{An example with $\omega \neq 0$} \label{nontrivialex}
Let us consider a non-holomorphic function $f: \D \ra \C \cup
\{\infty\} = \PP^1 \C$. Let $\gamma$ be the graph of the function
$f$ and then let us consider the manifold $X=(\D \times \C \PP^1)
\exc \gamma$. This manifold is provided with the natural complex
structure restricting the one of $\D \times \C \PP^1$ and it is
\hl ally foliated by the leaves $L_x=
\C \PP^1 \exc \{f(x)\} \simeq \C$.\\

Let $F: \D \times \C \ra X$ be the map in $\mathcal{C}$
constructed from a transverse \hle disk $D_0 \subset \D \times
\{0\}$, and from the \hle trivialization $t$ of the bundle
$T\mathcal{F}$ that is defined by $t_x=1 \, \in \C \simeq T_x \C
\PP^1 \simeq T_x T \mcl{F}$. Then, if we assume for example that
$f \neq 0$ and $f \neq \infty$ on $D_0$, $F$ reads as :
\begin{equation*}
F: \begin{cases} D_0 \times \C & \lgra \cup_{x \in D_0} L_x \\
                  (x,y)         &   \lgra \left(x, \f{f(x)
                  y}{y+f(x)}\right)
\end{cases}
\end{equation*}
\par It follows :
\begin{equation*}
dF_{(x,y)} \begin{pmatrix} v_x \\ v_y \end{pmatrix} = \begin{pmatrix}
v_x \\ \f{f(x)^2}{(y+f(x))^2} v_y + \f{y^2}{(y+f(x))^2}
\left(\f{\pp f}{\pp x}(x) v_x + \f{\pp f}{\pp \bar{x}}(x) \bar{v_x} \right) \end{pmatrix}
\end{equation*}

\par Thus, it is straightforward to check that $w(x,y)= 2i
\f{y^2}{f(x)^2} \f{\pp f}{\pp  \bar{x}}(x)$ and so does not
vanish, since  $f$ is not \hl.\\
\par Furthermore, let us note, that this is, in a way, ``the
easiest''
  example. Indeed, $\omega$ here is a polynomial of degree 2 with
  respect to $y$ ($\omega(x,y)=\f{2 i}{f(x)^2} \f{\pp f}{\pp
  \bar{x}}(x)\ y^2$), which can
  be considered as the ``easiest" example for a \hle function on $\C$
  that vanishes on 0 and whose derivative vanishes on 0.

\section{Final remarks and openings}
The Brody and Kobayashi-hyperbolicity, initially defined for
complex manifolds (see \cite{lang}), can also be defined in the
case of almost-complex structures. This issue has been tackled in
the last few years and is still an active area of research,
notably in the case of an almost-complex structure compatible with
a symplectic structure (see \cite{alb}, \cite{alb1}). However, the
notion of measure-hyperbolicity and the notion of hyperbolicity
for foliations defined in this paper only make sense in the
framework of integrable complex structures. Indeed, as soon as
there exists a \hle map from $\D^{n-1} \times \C$, or $\D^n$, to
an almost-complex manifold, then its structure has to be
integrable.

That's why, in the case of foliated manifolds provided with an
almost-complex structure (for example compatible with a \spe
structure), such that the leaves are complex submanifolds, it
would be interesting to tackle the issue of hyperbolicity
following different approaches. The one developed in \cite{alb}
and \cite{alb1} suggests us to use a notion of foliated Floer
homology.

\appendix
\section{Holomorphic foliations: definitions and notations}
\label{rappel} For more details, one can refer to \cite{god} and
\cite{molino}.\\

 A complex manifold $X$ is
provided with a {\it \hle foliation} $\mcl{F}$ if there exists an
atlas of \hle charts $\psi_i: U_i \ra \Omega_i \subset \C^p \times
\C^k$ such that the transitions maps are bi-holomorphisms:

\beq* \psi_i \circ \psi_j^{-1} :
\begin{cases} \psi_j^{-1} (U_i \cap U_j) & \lgra
\psi_i^{-1} (U_i \cap U_j) \\
(x,y) \in \C^p \times \C^k & \lgra (\phi^{i,j}_1(x),
\phi^{i,j}_2(x,y)) \in \C^p \times \C^k.
\end{cases}
\eeq*

Thus in each point, $\phi^{i,j}_1$ is a local diffeomorphism of
$\C^p$. And geometrically, if for each $x \in \C^p$ we call {\it
plaques} the connected components of $\psi_i^{-1}(\{x\} \times
\C^k)$, then the transition map sends each plaque (into $U_j$) to
another plaque (into $U_i$).

These charts are called {\it foliated charts}.

An open set $U$ on which there exists a foliated chart $\psi: U
\ra \D^p \times \D^k$ is called a {\it distinguished} open set.\\

In each point $z \in X$, the tangent plane to the foliation $(T
\mcl{F})_z$ is defined as $\dd \psi^{-1} (\{0\} \times \C^k)$ for
any foliated chart $\psi$. An equivalence relation can then be
defined: $x \sim y$ if $x$ and $y$ can be linked by a path tangent
to the plane field $T \mcl{F}$. The {\it leaves} of the foliation
are the equivalence classes for this relation.

They are the smallest connected sets satisfying: if a plaque
intersects this set, then this plaque is completely included into
it.

Furthermore, by definition, the foliated charts send each leaf to
a $\{*\} \times \C^k$.\\

The above introduced plane field $T\mcl{F}$ is a \hle vector
bundle above the manifold. Then, another \hle vector field can
naturally be defined: the { \it transverse or normal vector
field}: $ T \mcl{N}=TX / T\mcl{F}$.\\

Moreover, a subset of $X$ is said to be {\it saturated} if it is
the union of leaves. The {\it saturation} of an arbitrary subset
$A$ is the union of all leaves passing through $A$. Let us recall
that
the saturation of an open subset is an open set.\\

Furthermore, if $U$ is a distinguished open set, one can consider
the foliation induced by $\mcl{F}$ on $U$. This {\it induced
foliation} will be denoted by $\mcl{F}^U$.
The leaves of this foliation are the plaques of $U$.\\

Finally, let us introduce the projection along the leaves on a
transverse disk. Let $D$ be a transverse disk to the foliation. If
$U$ is a distinguished open set and $\psi: U \ra \D \times \D$ is
a foliated chart sending $D$ onto $\D \times \{0\}$, then we can
consider on $U$ the {\it projection $\pi$ on the transverse $D$
along the leaves}. This projection $\pi$ sends each plaque
$L^{U}_x$ (the leaf of the induced foliation on $U$ passing
through $\psi^{-1}(x,0)$) on $\psi^{-1}(x,0)$. Moreover, via  the
chart $\psi$, the map $\pi$ reads as \beq* \displaystyle \psi
\circ \pi \circ \psi^{-1}:
\begin{cases}
\D \times \D^p & \lgra \D^k \times \{0\} \\
z=(x,y) & \lgra z'=(x,0)
\end{cases} \eeq*
This is the projection on the first coordinate, so it is \hl.
Therefore, the projection $\pi$ on the transverse $D$
is \hle on $U$.\\

\section{Holomorphy of the projections along the leaves}
\label{projhol} Let $D$ be a fixed transverse \hle disk
parametrized by $\D$, and $F: \D \times \C \ra X$ a map in
$\mathcal{C}$ constructed along the disk $D$. Then
\begin{lemma}
\label{translationalongleaves} Let $D'$ be any transverse disk
included in the image of $F$, which can be parametrized as $D'=\{
F(x, \mcl{Y}(x)), x \in \D \}$, with $\mcl{Y}$ any continuous map
from $\D$ to $\C$. Then, on a distinguished neighborhood $U$ of
$D'$ can be defined a projection $\pi_0$ on $D$ satisfying: \beq*
\displaystyle \pi_0: \begin{cases} U & \lgra D \\
z=F(x, \mcl{Y}(x)) & \lgra F(x,0),
\end{cases} \eeq*
It also sends the plaque of $U$ passing through $z=F(x,
\mcl{Y}(x))$ onto $F(x,0)$.\\
Moreover, it is \hle on $U$.
\end{lemma}

\begin{proof}
On any distinguished open set $U_0$, one can consider the
projection $\pi_0$ on $D$ along the leaves, defined the above
section. This projection is \hle and sends each plaque of $U_0$
passing through $F(x,0)$ on $F(x,0)$. If $D'$ is included in $U_0$
(with a ``small" parametrization $\mcl{Y}$) this projection
satisfies the wished properties.

\noindent In order to deal with a general $D'$, the idea (used
also to define the holonomy pseudo-groups) is to link $D'$ with
$D$ (along for example $\left(F(x, t \mcl{Y}(x)\right)_{t \in
[0,1]}$) through a chain of distinguished open sets $(U_j)_{j=0
\ldots n}$, $D \subset U_0$, $D' \subset U_n=V$, such that for any
$j= 1 \ldots n$, $U_{j-1} \cap U_j$ contains a transverse \hle
disk $D_j$.\\

\noindent Then on each $U_j$, we can, as above, define a \hle
projection $\pi_j$ on the disk $D_j$. Then we can define on
$V=U_n$ the projection $\pi_0 \circ \pi_1 \circ \ldots \circ
\pi_n$. It is \hle and satisfies the desired properties.
\end{proof}

Let us notice that this projection map on $D$ is not necessarily
unique since it might depend on the choice of the path linking
$D'$ and $D$.\\

The idea of this lemma is similar to the one used to define the
holonomy groups of the leaves: if $L_0$ is a leaf, and $\gamma$ is
an element of $\pi_1(L_0)$ represented by a loop $\gamma_0 \subset
L_0$ with base-point $x_0$, then one can consider a chain of
distinguished open sets $(U_j)$ along $\gamma_0$, with some
transverse \hle disks $D_j$ in the intersection $U_j \cap
U_{j-1}$. Then defining $\pi_j$ the projection on $D_j$ along the
leaves -- defined on $U_j$ --, one can consider, as above,
$\pi_\gamma$ the composition of these projections: $\pi_\gamma :
D_0 \ra D_0$ with $D_0$ a small transverse \hle disk around $x_0$.
As seen above, this map is holomorphic (and even a local
bi-holomorphism around $x_0$ of inverse $\pi_{-\gamma}$). This
determines a germ of bi-holomorphism around $x_0$ which does not
depend on the choices of the chain $(U_j)$ and of the
representation $\gamma_0$. This defines a map from $\pi_1(L_0,
x_0) \ra Hol(D_0)$ with $Hol(D_0)$ the group of germs of local
bi-holomorphism around $x_0$. The image of this map is the
holonomy group of $L_0$ at $x_0$. For two distinct points in
$L_0$, the holonomy groups of $L_0$ in these points are
conjugated. So one can talk about the
holonomy group of the leaf $L_0$.\\

Let us come back to the result of lemma
\ref{translationalongleaves}. From this lemma it immediately
follows the existence and holomorphy of the function $\theta_1$
introduced in section \ref{changt}.
\par
Moreover, we can deduce:

\begin{lemma} The complex
structures $J'=F^{*}J$ and  $J_0$ of $\D \times \C$ coincide on
the quotient $\D$. \end{lemma}

\begin{proof}
On a (distinguished) neighborhood $V$ of a fixed point $z=(x_0,
y_0)$, the projection $\pi_0$ (introduced in the previous
proposition) can be defined. It satisfies $\pi_0(F(x,y))=F(x,0)$
on $V$, and so differentiating

\be \label{etunepetitedeplus} \dd(\pi_0)_{F(x,y)}(\dd
F_{(x,y)}(v,0))=\dd F_{(x,0)}(v,0).\ee Since $\pi_0$ is \hl,
$\dd(\pi_0)_{F(x,y)}(J Y)= J \dd(\pi_0)_{F(x,y)}(Y)$. Applying
this equality to $\dd F_{(x,y)}(v,0)$, we get
  $\dd(\pi_0)_{F(x,y)}(J\dd F_{(x,y)}(v,0))= J \dd(\pi_0)_{F(x,y)}(\dd
F_{(x,y)}(v,0))$, which equals $J \dd F_{(x,0)}(v,0)$ according to
(\ref{etunepetitedeplus}). Considering that $F$ is \hle along
     $\D \times \{ 0\}$, this  equals $\dd F_{(x,0)}(J_0 v,0)$.\\
Finally, applying (\ref{etunepetitedeplus}) one last time, \beq*
\dd(\pi_0)_{F(x,y)}(J\dd F_{(x,y)}(v,0))=
 \dd(\pi_0)_{F(x,y)}(\dd F_{(x,y)}(J_0 v,0)). \eeq* Since
$\mathrm{Ker}(\dd\pi_0)=T\mathcal{F}$, our result follows.
\end{proof}


\begin{thebibliography}{99}


\bibitem{alb}
A.L. Biolley, {\textit{Cohomologie de Floer, hyperbolicit\'es \spe
et \pc}}, Th{\`e}se de Doctorat de l'Ecole Polytechnique (2003).

\bibitem{alb1} A.L. Biolley, {\it Floer homology,
symplectic and complex hyperbolicity}, math.SG/0404551, to be
published (2004).

\bibitem{edsull}
R. Edwards, K. Millett, D. Sullivan, {\it Foliations with all
leaves compact}, Topology {\bf 16} (1977).


\bibitem{ghys}
E. Ghys,  {\it Dynamique et g{\'e}om{\'e}trie complexes},
Panoramas et Synth{\`e}ses {\bf 8} (1997).


\bibitem{god} C. Godbillon, {\it Feuilletages, Etudes
g\'eom\'etriques}, Progress in Mathematics, Birkh\"auser (1991).

\bibitem{holmann} H. Holmann, {\it On the stability of holomorphic
foliations with all leaves compact},  Variétés analytiques
compactes (Colloq. Nice, 1977), Lecture Notes in Mathematics {\bf
683}, Springer, (1978).


\bibitem{holmann2} H. Holmann, {\it On the stability of
holomorphic foliations},  Analytic functions, Kozubnik 1979 (Proc.
Seventh Conf., Kozubnik, 1979), Lecture Notes in Mathematics {\bf
798}, Springer, (1980).



\bibitem{lang} S. Lang, {\it Introduction to complx hyperbolic
spaces},Springer-Verlag (1987).


\bibitem{molino} P. Molino, {\it Riemannian Foliations}, Progress
in Mathematics, Birkh\"auser (1988).


\bibitem{nishino} T. Nishino, {\it Nouvelles recherches sur les fonctions enti\`eres
    de plusieurs variables complexes (II). Fonctions enti\`eres qui se
    reduisent \`a celles d'une variable},  J. Math. Kyoto Univ
    {\bf 9-2} (1969).

 \bibitem{il1} Yu.S. Ilyashenko, {\it Foliations by analytic
 curves}, Mat. Sb. {\bf 88(130)} (1972).

 \bibitem{il2} Yu.S. Ilyashenko, {\it Covering manifolds for analytic
 families of leaves of foliations by analytic curves}, Topol. Methods
 Nonlinear Anal. {\bf 11, no2} (1998).

\end{thebibliography}
\end{document}